\documentclass[%
 aip, 
cp,  
 amsmath,amssymb,
 reprint,%
]{revtex4-2}

\usepackage{graphicx}
\usepackage{dcolumn}
\usepackage{bm}

\usepackage[utf8]{inputenc}
\usepackage[T1]{fontenc}
\usepackage{mathptmx} 

\begin{document}

\title{Time-dependent finite-dimensional dynamical system representation of breather solutions} 

\author{Yoritaka Iwata} 
\email[Corresponding author: ]{iwata_phys@08.alumni.u-tokyo.ac.jp}
\affiliation{%
Faculty of International Studies, Osaka University of Economics and Law 
}

\author{Yasuhiro Takei}

\affiliation{Mizuho Research \& Technologies}

\date{\today} 

\begin{abstract}
A concept of finite-dimensional dynamical system representation is introduced.
Since the solution trajectory of partial differential equations are usually represented within infinite-dimensional dynamical systems, the proposed finite-dimensional representation provides decomposed snapshots of time evolution.
Here we focus on analyzing the breather solutions of nonlinear Klein-Gordon equations, and such a solution is shown to form a geometrical object within finite-dimensional dynamical systems.
In this paper, based on high-precision numerical scheme, we represent the breather solutions of the nonlinear Klein-Gordon equation as the time evolving trajectory on a finite-dimensional dynamical system.
Consequently, with respect to the evolution of finite-dimensional dynamical systems, we confirm that the rotational motion around multiple fixed points plays a role in realizing the breather solutions.
Also, such a specific feature of breather solution provides us to understand mathematical mechanism of realizing the coexistence of positive and negative parts in nonlinear systems.
\end{abstract}

\maketitle

\section{Introduction}
Due to the underlying nonlinearity, natural events possibly take place in a complex way.
For example, pattern formation, which can be found in physics, chemistry, biology, and materials science, typically arises from the nonlinearity. 
Also nonlinear wave propagation such as shock waves and solitons are observed in different scales~\cite{93Cross, 20Nabika, 15Ball}.
In order to understand these nonlinear phenomena quantitatively, we introduce a mathematical tool to represent and analyse the time evolution driven by nonlinearity.

In this paper, we propose a method to represent a point (the state of the system at a certain time) on a solution trajectory inside the infinite-dimensional dynamical system as a geometrical object in finite-dimensional dynamical system.
We call this method "$t$-dependent finite-dimensional representation" of originally infinite-dimensional dynamical system (for a textbook of finite-dimensional dynamical systems see \cite{74hirsch}, and for a textbook of infinite-dimensional dynamical systems see \cite{97temam}).
Here note that solution trajectories of partial differential equations are usually represented within infinite-dimensional dynamical systems.
The method is then applied to clarify the characteristics of breather mode of the nonlinear wave described by nonlinear Klein-Gordon equations.
As a result, the trajectory of breather solution provides a coexistence of positive and negative parts, which are geometrically represented in finite-dimensional dynamical systems by coexistence of rotational motions around multiple fixed points.

\section{Finite-dimensional representation of infinite-dimensional dynamical systems}
\subsection{Infinite-dimensional dynamical system}
A couple of space and semigroup $(X, S(t))$ is called dynamical system.
Here we limit ourselves to the solvable cases with a well-defined semigroup $S(t)$ on a Banach space $X$, and concentrate on introducing a tool analyzing nonlinear dynamics in a geometrical manner.
One parameter semigroup of operator satisfies the semigroup property
\[ \left\{ \begin{array}{ll}
S(t+s)  = S(t) ~ S(s), \quad t,s \ge 0, \vspace{3mm} \\
S(0) = I.
\end{array}  \right. \]
We consider a Cauchy problem of nonlinear wave equations.
It is formulated as abstract evolution equations of hyperbolic type
\begin{equation} \label{eq01}
\begin{array}{ll} 
d U(t)/dt  +  
\left(
\begin{array}{ll}
~ 0  \qquad  -I \\
-\partial_x^2 \qquad  0
\end{array}
\right)
U(t) = 
\left(
\begin{array}{cc}
0  \\
-\mu^2 u(t) + k(u(t)) 
\end{array} 
\right),  \vspace{5mm} \\
U(0) = U_0 = (u_0, v_0)^T,
\end{array}
\end{equation}
in an infinite-dimensional Banach space $X$, where the unknown function $U(t)$ is represented by a vector form $U(t) = (u(t), v(t))^T$ with a relation $v(t) = \partial_t u(t)$.
For a real positive number $\mu$ standing for a mass, the non-autonomous term $F(U(t))$ consists of linear part $(0,-\mu^2 u(t))^T$ and the nonlinear part $(0, k(u(t)))^T$, where $k(u(t))$ is, for example, a polynomial of $u(t)$ and/or a trigonometric function of $u(t)$.
Let the space $X$ consist of square-integrable periodic functions being defined on one-dimensional interval; $X = L_{\rm per}^2(0,L) \times  L_{\rm per}^2(0,L)$, and semigroup $S(t): ~ U(s)\mapsto U(t+s)$ with $t,s \ge 0$ defined on $X$.
If it comes to partial differential equations, the dynamical system $(X, S(t))$ becomes infinite-dimensional, as $X = L_{\rm per}^2(0,L) \times  L_{\rm per}^2(0,L)$ is infinite-dimensional.


\subsection{Finite-dimensional dynamical system}
The finite-dimensional representation is introduced with a focus on the spatial variables.
Such a viewpoint naturally appears, if we limit ourselves to the stationary problem of originally non-stationary problem.
The stationary problem of (\ref{eq01}) is given by
\begin{equation} \label{eq02} \begin{array}{ll}
\begin{array}{ll} 
\left(
\begin{array}{ll}
~ 0  \qquad  -I \\
-\partial_x^2 \qquad  0
\end{array}
\right)
\end{array}
 U(t) = 
\left(
\begin{array}{cc}
0  \\
-\mu^2 u(t) + k(u(t)) 
\end{array} 
\right)
\end{array}
\end{equation}
in $X$ at a certain fixed $t$, which is equivalent to 
\begin{equation} \label{eq03}
\left\{
\begin{array}{ll}
 \partial_x^2 u(t,x) = \mu^{2} u(t,x) - k(u(t,x)),  \vspace{3mm}  \\
 u(t,0)  =  u(t,L), \quad   \partial_x u(t,0)  =  \partial_x u(t,L),  
\end{array}
\right.
\end{equation} 
in $(0,L)$.
Equation (\ref{eq03}) is an ordinary differential equation holding a finite-dimensional dynamical system whose evolution direction is intentionally taken to be the spatial direction.
Here it is practical to introduce two parameter semigroup
\[
s(x_2,x_1): ~ u(x_1) ~ \mapsto ~ u(x_2),
\]
on Euclidean space ${\mathbb R}^2$.
Two-parameter semigroup $s(x_2,x_1)$ satisfies the semigroup property
\[ \left\{ \begin{array}{ll}
s(x_3,x_1)  = s(x_3,x_2) ~ s(x_2,x_1), \quad 0 \le x_1 \le x_2 \le x_3 \le L,  \vspace{3mm} \\
s(x,x) = I.
\end{array}  \right. \]
Since the periodic boundary condition is assumed to be satisfied in the present settings, the finite-dimensional dynamical system shows a closed loop, where note that ``closed" loop appears due to the periodic boundary condition.
The loop is not necessarily a simple closed curve, and some crossings possibly appear.
According to the possible appearance of crossing points, here we introduce a two-parameter semigroup $s(x_2,x_1)$, in which the evolution can be different depending on both the starting point $x_1$ and ending point $x_2$.

The finite-dimensional representation in this paper is obtained based on the stationary problem; $\partial_t u = 0$, so that it is regarded as a snapshot at a certain time.
That is, the solution at a certain time is shown in the $u-v$ plane (i.e., phase space).
If those snapshots are connected with respect to the variable $t$, then a $t$-dependent moving loop is obtained.
It is called here the time-dependent finite-dimensional representation.
Although the finite-dimensional representation is proposed here particularly for nonlinear hyperbolic evolution equations, it can be applied to any partial differential equations including fluid dynamics.

\section{Visualization of finite-dimensional dynamical systems}
High-precision numerical scheme is indispensable in the numerical treatment of nonlinear problems.
In order to obtain the reliable time-evolving geometric shapes, we employ high-precision numerical scheme satisfying the hyperbolic type conservation laws; for benchmark tests, see~\cite{21takeiiwata-2,20iwata-2,21takeiiwata-3}.
In the high-precision numerical scheme, the spatial discretization is made by the Fourier spectral method (for example, see \cite{06canuto,97Gottlieb,86Canuto}), and the time-discretization by the implicit Runge-Kutta method (for example, see \cite{87Butcher,20Rana,07Mardal,06Staff}).

\begin{figure*}[t]
\begin{center}
  \includegraphics[width=52mm,bb=9 9 358 234]{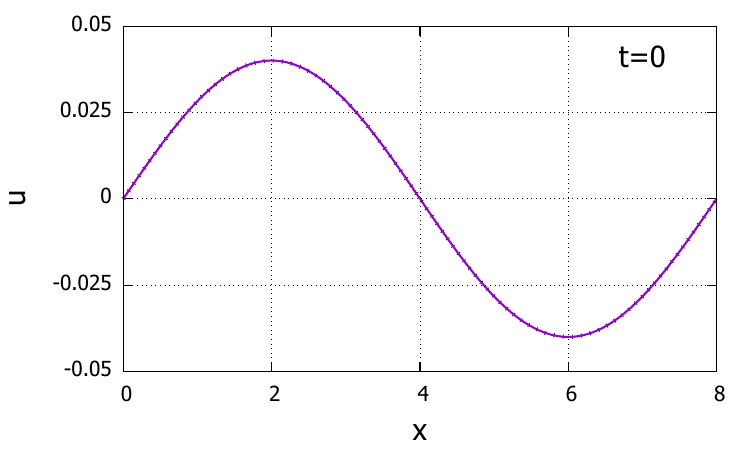} 
  \includegraphics[width=52mm,bb=9 9 358 234]{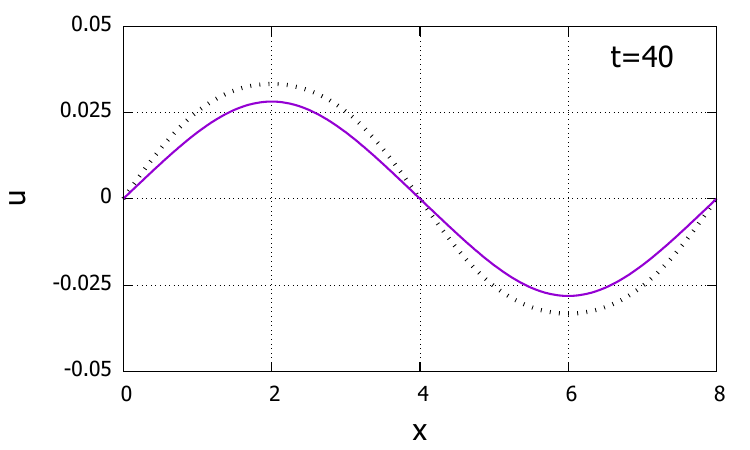} 
  \includegraphics[width=52mm,bb=9 9 358 234]{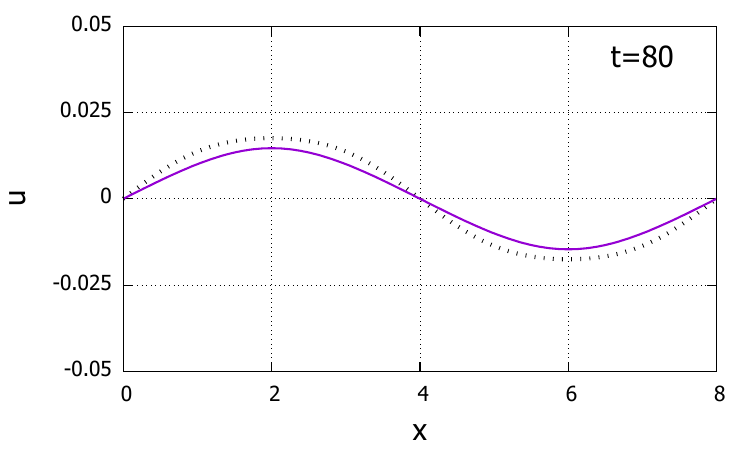}  \\
  \includegraphics[width=52mm,bb=9 9 358 234]{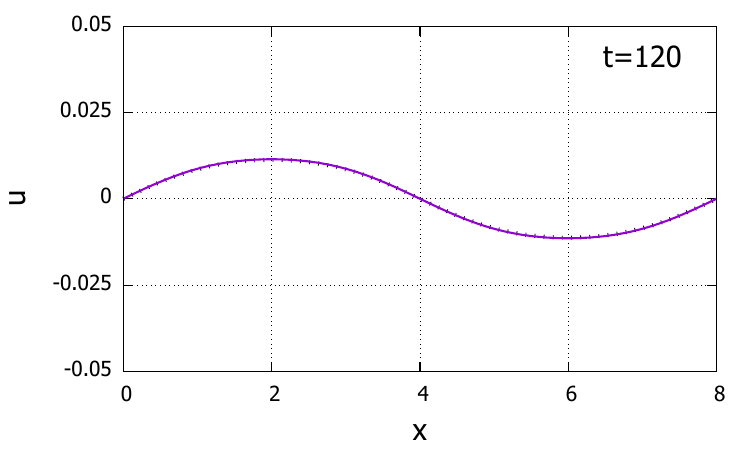} 
  \includegraphics[width=52mm,bb=9 9 358 234]{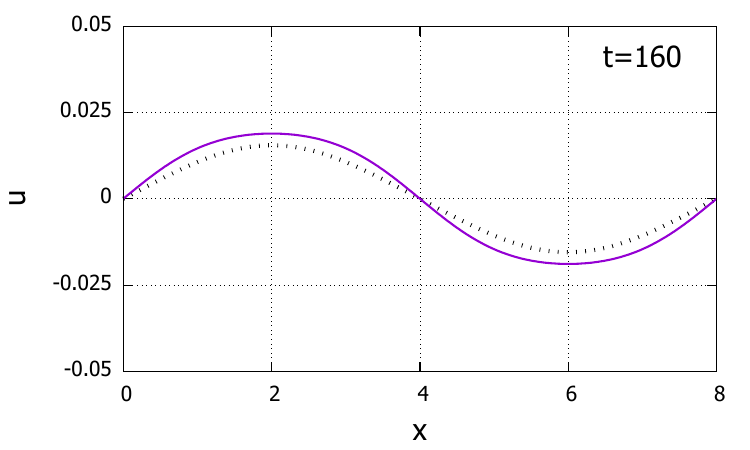} 
  \includegraphics[width=52mm,bb=9 9 358 234]{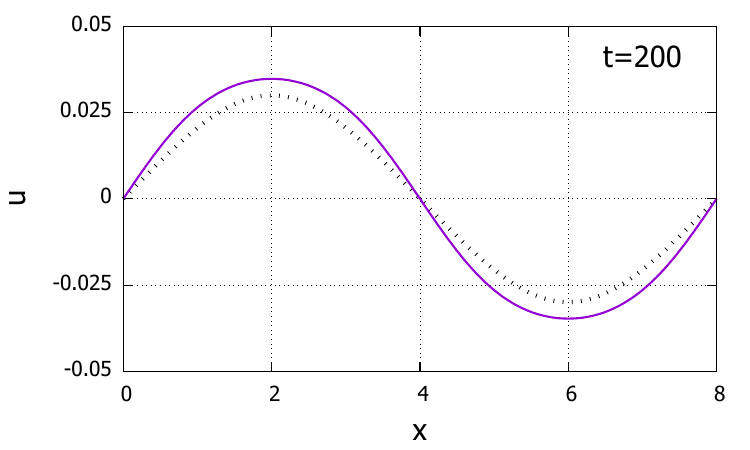} 
\caption{Breather solution,\ $(A = 0.04)$ : \ 
In each panel, the solid line shows the function $u$ at each time $t$, and the dotted line shows the function $u$ at time $(t-16)$.
A part $u(t,x)$ with $x \le 4$ is confined in the positive region ($u(t,x)>0$) and a part $u(t,x)$ with $x > 4$ is in the negative region ($u(t,x)<0$).
}
\label{fig-b}
\end{center}
\end{figure*}

\subsection{Nonlinear Klein-Gordon equations with cubic nonlinearity}
Let $x \in [0,L]$ be a finite domain of space.
The time evolution problem is considered for $t \ge 0$. 
We consider the nonlinear Klein-Gordon equation with the double-well type interaction.
\[
\label{eq:eq301}
\begin{array}{ll}
\ \partial_t^{2} u + \alpha \partial_x^{2} u +   (\beta u^2 - \mu) u  = 0, \vspace{3mm} \\
\ u(0,x) = A \sin(\pi x/4),\quad u(t,0) = u(t,L),  \vspace{3mm} \\
\ \partial_t u(0, x) = 0,\ \partial_t u(t, 0) = \partial_t u(t, L).
\end{array}
\qquad {\rm (KG)}  
\]
where $\alpha = 2^{-8}$, $\beta = 1$, $\mu = 0.00305$, and $A>0$ are real constants.
This equation is known to hold the symmetry breaking being known as the Higgs mechanism; $\phi^4$-theory in the context of quantum field theory (for a textbook, see \cite{65bjorken}).

\subsection{Ordinary oscillatory solution and the breather solution}
Two different solutions (sometimes a terminology ``modes" is used for ``solutions") appears for the Cauchy problem (KG) only by changing an amplitude $A$ of the initial function. 
One is the ordinary oscillatory mode and the other is the breather mode, where it is practical to imagine harmonic oscillation for the ordinary oscillatory solution.
The ordinary oscillatory mode and the breather mode are the periodic solutions for both $t$ and $x$, so that they are a kind of doubly periodic and therefore closed curves as for both space and time.

The breather solution can be distinguished from an ordinary oscillatory solution by the appearance of certain kinds of collectivity, where activated modes result in the resonance.
For the collectivity, the breather solution includes the resonating large amplitude oscillation, which is localized only in the positive or negative side of $u=0$ (Fig.~\ref{fig-b}).
Such a localization is not seen in the ordinary oscillatory solution.
The breather mode is achieved by the collective resonance of many activated modes, as such curved cannot be made of a sum of few sine or cosine curves.

\begin{figure*}[t]
\begin{center}
  \includegraphics[width=40mm,bb=65 20 295 200]{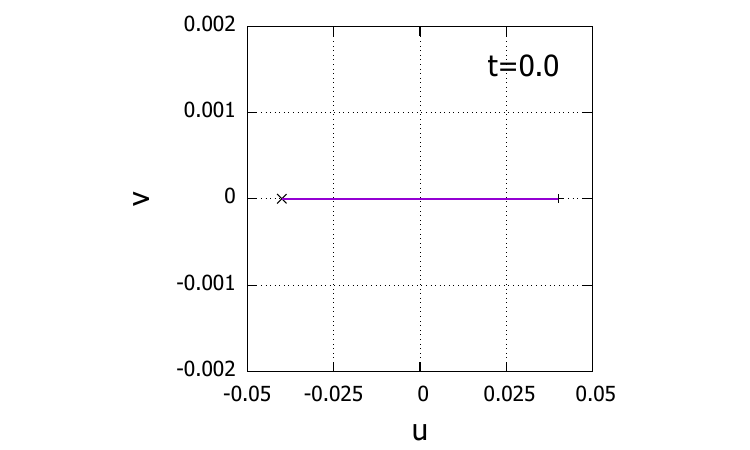} 
  \includegraphics[width=40mm,bb=65 20 295 200]{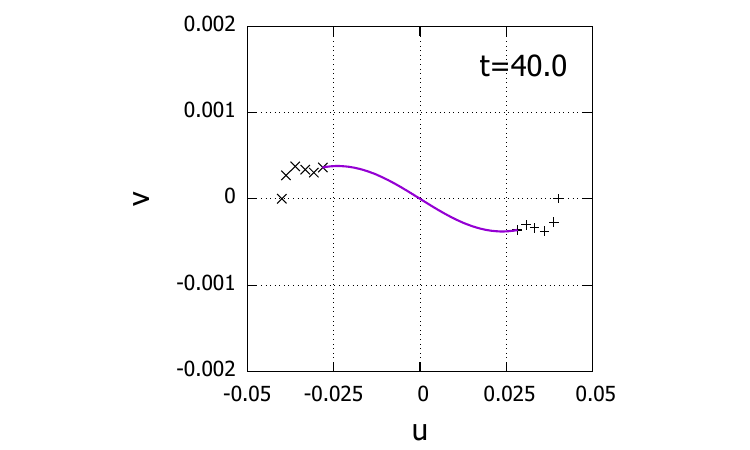} 
  \includegraphics[width=40mm,bb=65 20 295 200]{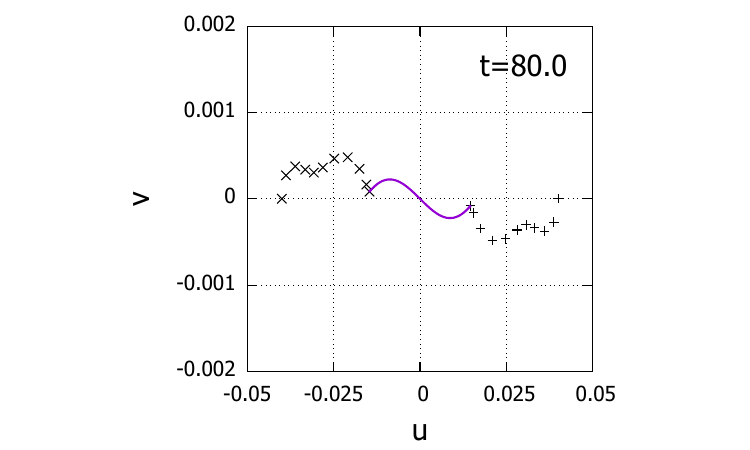}  \vspace{8mm} \\
  \includegraphics[width=40mm,bb=65 20 295 200]{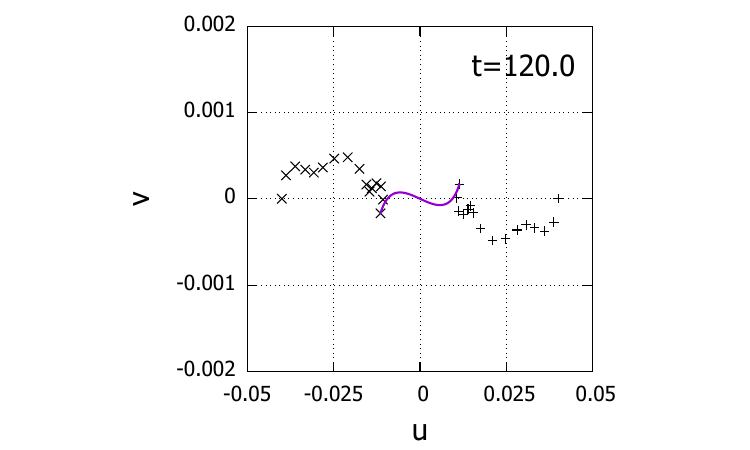} 
  \includegraphics[width=40mm,bb=65 20 295 200]{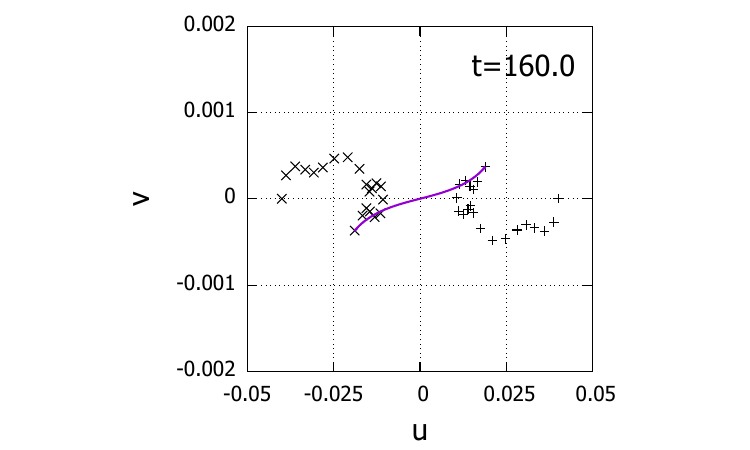} 
  \includegraphics[width=40mm,bb=65 20 295 200]{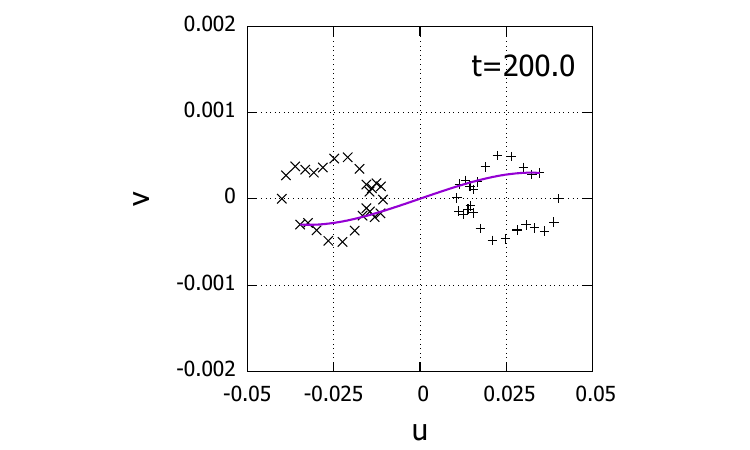}  \vspace{4mm}
\caption{The finite-dimensional representation for a breather solution $(A = 0.04)$ : \ 
The solid line represents the finite-dimensional representation for the breather solution at each time $t$, represented by a curve on the $u-v$ plane.
The points on the graph indicate the trajectories of the endpoints of the curve up to time $t$.
The left and right halves of the curve evolve in time to rotational motion with the fixed point $( \pm (\mu/\beta)^{1/2}, 0)$ as the focus, respectively.
}
\label{fig-b_rep}
\end{center}
\end{figure*}



Figure \ref{fig-b_rep} shows a $t$-dependent finite-dimensional representation of the breather solution.
At $t=0$, the solution is represented as a line segment along the $u$-axis on the $u-v$ plane, as can be seen from the initial value setting.
For the ordinary oscillatory solution, the curve rotates around the origin $(0,0)$.
Unlike the ordinary oscillatory solution, the curve does not rotate around the origin $(0,0)$ in case of breather solution, and furthermore, the length of segments becomes shorter and larger.
Instead additional two fixed point $(\pm (\mu/\beta)^{1/2},0)$ play roles.

\section{Summary}
We propose $t$-dependent finite-dimensional representation of originally infinite-dimensional dynamical systems.
The proposed method is applied to a breather solution for nonlinear Klein-Gordon evolution equation with cubic nonlinearity.
By comparing Figs.~\ref{fig-b} and \ref{fig-b_rep}, more details of breather solution can be seen by having a finite-dimensional representation (Fig.~\ref{fig-b_rep}). 
Although ordinary oscillatory solution is a rotation centered at the origin $(0,0)$ on the $u-v$ plane, the breather solution is more complex. 
Indeed, for the breather solution, two local rotational motions coexisting at the fixed point $(\pm (\mu/\beta)^{1/2},0)$.

\nocite{*}

\end{document}